\documentclass[12pt]{article}

\usepackage{latexsym}
\usepackage{amssymb}
\usepackage{graphicx}
\usepackage{color}
\usepackage{amsmath}
\usepackage{mathrsfs} %fonts veramente calligrafici con \mathscr{B}

\newtheorem{theorem}{Theorem}

\newtheorem{lemma}[theorem]{Lemma}

\newcommand{\kom}[1]{}
\renewcommand{\kom}[1]{{\bf [#1]}}
\definecolor{blau}{rgb}{0.1,0.0,0.9}

\newcounter{komcounter}
\numberwithin{komcounter}{section}

\begin{document}

\title{Growth of subsolutions to fully nonlinear equations in halfspaces II}

\author{FATMA GAMZE D\"UZG\"UN{\small{$^1$}}\\ NIKLAS L.P. LUNDSTR\"OM{\small{$^2$}} and MARCUS OLOFSSON{\small{$^3$}}
%Niklas L. P. Lundstr\"{o}m %{\small{$^2$}}
\\\\
{\small{$^1$}}
\it \small Department of Mathematics and Computer Science, \\ 
\it \small University of Cagliari, Cagliari, Italy.\\
{\small{$^2$}}
\it \small Department of Mathematics and Mathematical Statistics,\\
\it \small Ume{\aa} University, SE-90187 Ume{\aa}, Sweden\/{\rm ;}\\
\it \small niklas.lundstrom@umu.se;\\
{\small{$^3$}}
{\it \small Odensalaskolan, Östersund, Sweden.
%bakkenolofsson@gmail.com
}  
}

\maketitle

%%%%%%%%%%%%%%%%%%%%%%%%%%%%%%%%%%%%%%%%%%%%
%%%%%%%%%%%%%%%%%%%%%%%%%%%%%%%%%%%%%%%%%%%%
%%%%%%%%%%%%%%%%%%%%%%%%%%%%%%%%%%%%%%%%%%%%
%%%%%%%%%%%%%%%%%%%%%%%%%%%%%%%%%%%%%%%%%%%%
%%%%%%%%%%%%%%%%%%%%%%%%%%%%%%%%%%%%%%%%%%%%
%%%%%%%%%%%%%%%%%%%%%%%%%%%%%%%%%%%%%%%%%%%%
%%%%%%%%%%%%%%%%%%%%%%%%%%%%%%%%%%%%%%%%%%%%
%%%%%%%%%%%%%%%%%%%%%%%%%%%%%%%%%%%%%%%%%%%%
%%%%%%%%%%%%%%%%%%%%%%%%%%%%%%%%%%%%%%%%%%%%
%%%%%%%%%%%%%%%%%%%%%%%%%%%%%%%%%%%%%%%%%%%%

\begin{abstract}
We establish sharp lower growth estimates for subsolutions
of fully nonlinear elliptic equations in half-spaces,
allowing source terms that depend explicitly on the solution.
The estimates capture both growth at infinity and blow-up
at a finite distance from the boundary. Our approach extends
the ODE-based method introduced in
[J. Differential Equations 320 (2022), 143--173].\\

\noindent
{\em Mathematics Subject Classification:}
% this was used on the JDE-paper:
35B40, 35B50, 35B53, 35D40, 35J25, 35J60, 35J70.

\end{abstract}

\section{Introduction and statement of main results}

\setcounter{theorem}{0}
\setcounter{equation}{0}

%KPP–Fisher equation  

%https://www.sciencedirect.com/topics/earth-and-planetary-sciences/reaction-diffusion-equation

%A nonlinear diffusion problem with convection and anisotropic nonstandard growth conditions

\noindent
In \cite{L22}, the second author introduced an approach to deriving sharp Phragm\'en--Lindel\"of-type 
theorems for fully nonlinear elliptic equations in half-spaces, reducing lower growth estimates for subsolutions to the analysis of associated ODEs. 
Here we extend this approach to more general lower-order
terms depending on both the solution and its gradient, and
obtain estimates covering both growth at infinity and blow-up
at a finite distance from the boundary.

We consider fully nonlinear nonhomogeneous degenerate elliptic partial differential equations (PDEs) in nondivergence form,
\begin{align}\label{eq:main}
F(x,u,Du,D^2u) = 0,\tag{$\star$}
\end{align}
in halfspaces of $\mathbb{R}^n$, $n \geq 1$.
Here, $Du$ is the gradient, $D^2u$ the Hessian,
$F:\mathbb{R}^n \times \mathbb{R} \times \mathbb{R}^n  \times \mathbb{S}^n \rightarrow \mathbb{R}$ in which $\mathbb{S}^n$ is the set of symmetric $n \times n$ matrices equipped with the positive semi-definite ordering;
for $X, Y \in \mathbb{S}^n$, we write $X \leq Y$ if $\langle (X - Y) \xi, \xi \rangle \leq 0$ for all $\xi \in \mathbb{R}^n$.
Without loss of generality we fix the halfspace to $\mathbb{R}^n_+ := \{x \in \mathbb{R}^n : x_n > 0\}$ and for a point $x \in \mathbb{R}^n$ we use the notation $x = (x_1, x_2, \dots, x_{n-1}, x_n) = (x', x_n)$. 
The operator $F$ is assumed to satisfy
\begin{align} \label{eq:ass_drift_super}
F(x,u,p,X) &\geq  F(x,v,p,Y) \quad \text{whenever} \quad u \geq v, X \leq Y,  \quad \text{and}  \notag\\
-F(x,r,p,X) &\leq \Phi(t ,r,|p|) + \Lambda(x_n) \text{Tr}(X^+) - \lambda(x_n) \text{Tr}(X^-) \tag{$\star\:\star_{\:F}$} 
\end{align}
whenever $r \in [0,\infty), x, p \in \mathbb{R}^n, X \in \mathbb{S}^n$, $X = X^+ - X^-$, $X^+ \geq 0$, $X^- \geq 0$ and $X^+  X^- = 0$.
Here, \(\lambda, \Lambda:[0,\infty) \to (0,\infty)\) are continuous, 
$\lambda$ nonincreasing, $\Lambda$ nondecreasing, and $0<\lambda\leq\Lambda<\infty$.
Concerning $\Phi$ and $t$ we assume that \(\Phi:[0,\infty)\times(-\infty,\infty)\times[0,\infty)\to\mathbb R\) is continuous and that either
\begin{align}\label{eq:ass_phi}
&\text{$t=|x|$ and $\Phi$ is nonincreasing in its first argument, or}\notag\\
&\text{$t=x_n$ and $\Phi$ is nondecreasing in its first argument.}\tag{$\star\:\star_{\:\Phi}$} \\
&\text{In either case, it holds that $\Phi\geq0$
or $\Phi\leq0$ on
$[0,\infty)\times[0,\infty)\times[0,\infty)$.}\notag
\end{align}
Note that the above assumptions allow the ellipticity to blow up as $\lambda$ may vanish, and $\Lambda$ may explode, at infinity.

When $\Phi \geq 0$, we will make use of solutions of the following second order initial value problem 
\begin{align}\label{eq:ODE}
\frac{df}{dt} &= -  \frac{\Phi(t,h(t),f(t))}{\lambda(t)} - K(R) \frac{\Lambda(t)}{\lambda(t)} f(t), \tag{$\star\star\star$} \\  \notag
\frac{dh}{dt} &= f(t),
\quad t \in(0,R), \quad h(0)= 0,\quad f(0) = \nu,
\end{align}
for some $\nu > 0$ and $R > 0$. 
If $\Phi\leq 0$ we use instead solutions of \eqref{eq:ODE} with $\lambda(t)$ replaced by $\Lambda(t)$ in the first term on the right hand side. 

Throughout the paper, we denote by $h_{\nu,R} = h_{\nu,R}(t)$ and $f_{\nu,R} = f_{\nu,R}(t)$ a solution of \eqref{eq:ODE} with $K(R) = n/{\gamma(R)}$ and 
 $h_{\nu} = h_{\nu}(t)$ and $f_{\nu} = f_{\nu}(t)$ a solution of \eqref{eq:ODE} with $K(R) \equiv 0$.
We will also make use of a positive nondecreasing function $\gamma = \gamma(R) : (0,R_\infty)  \rightarrow \mathbb{R}_+$, which we will choose depending on the situation, determining the necessary properties of the domain
\begin{align}\label{eq:domain_2}
D(R) := B((0,\dots,0,-\gamma), \gamma + R) \cap \mathbb{R}_+^n,
\end{align}
where $B(x,r)=\{y\in\mathbb{R}^n:|x-y|<r\}$.
Moreover, for a subsolution $u$ and $R > 0$ we define  
\begin{align*}
M(R) = \sup_{\partial D(R)} u \quad\text{and}\quad M'(R) = \liminf_{\delta \to 0^+} \frac{M(R) - M(R-\delta)}{\delta}.
\end{align*}
The following theorem, which is a generalization of Lundstr\"om \cite[Theorem 2.1]{L22}, gives a lower growth estimate of subsolutions of
\eqref{eq:main} in terms of the solutions $f_{\nu,R}$ and $f_{\nu}$ of the ODE \eqref{eq:ODE}:
\begin{theorem}
\label{th:Phragmen--Lindelof}
Suppose that \eqref{eq:ass_drift_super} and \eqref{eq:ass_phi} hold and let $u$ be a subsolution of \eqref{eq:main} in $\mathbb{R}^n_+$ satisfying
\begin{align*}%\label{eq:tjohejsats}
\limsup_{x \to y} \, u(x) \, \leq \,  0 \quad \textrm{for all} \quad y \in \partial \mathbb{R}^n_+. %\quad \textrm{and} \quad M(1) \geq u_0,
\end{align*}
Let $\nu>0$ and $R_\infty\in(0,\infty]$.
Assume that the initial value problem \eqref{eq:ODE}
with $K\equiv 0$ has a solution satisfying $f_\nu>0$
on $(0,R_\infty)$ and that, for every
$R\in(0,R_\infty)$, the initial value
problem with $K=K(R)$ has a solution on $[0,R]$ satisfying 
$f_{\nu,R}>0$ and $h_\nu>h_{\nu,R}$
on $(0,R)$.
If $u(\bar x) \geq  h_{\nu}(\bar{x}_n)$ for some $\bar x = (0,\bar x_n)$ with $0 < \bar x_n <R_\infty$, 
then $M(R)$ is increasing on $(\bar x_n, R_\infty)$ and
\begin{align*}%\label{eq:limit main}
\liminf_{R\to R_\infty} \, \frac{M'(R)}{f_{\nu}(R)} \, \geq \,\liminf_{R\to R_\infty} \, \frac{f_{\nu, R}(R)}{f_{\nu}(R)}.
\end{align*}
\end{theorem}
We remark that if we can prove
$$
\liminf_{R\to R_\infty} \, \frac{M'(R)}{f_{\nu}(R)}  > 0,
$$
then
$
M(R) - M(R_0) \geq c \int_{R_0}^{R} h'_{\nu}(t) dt = c(h_{\nu}(R) - h_{\nu}(R_0))
$
whenever $R > R_0$ for some positive $R_0$, $c$ and thus %it follows that
\begin{align}\label{eq:hejdutydligagrowthest}
\liminf_{R\to R_\infty} \, \frac{M(R)}{h_{\nu}(R)} > 0.
\end{align}
Hence, in such case, Theorem \ref{th:Phragmen--Lindelof} takes the form of classical Phragm\'en--Lindel\"of theorems.
We also remark that Theorem \ref{th:Phragmen--Lindelof} 
also applies to domains $\Omega\subset\mathbb R^n_+$, 
provided that the boundary condition is imposed on all of $\partial\Omega$ and
$M(R)$ is defined accordingly.

The growth estimate in Theorem \ref{th:Phragmen--Lindelof} is sharp 
whenever the associated 
one-dimensional solution $u(x) = h_\nu(x_n)$ 
gives an admissible example. 
Indeed, consider the model operator
\begin{align}\label{eq:sharp}
 F(x,r,p,X) = -\Lambda(x_n)\operatorname{Tr}(X^+) + \lambda(x_n)\operatorname{Tr}(X^-) - \Phi(x_n,r,|p|),
\end{align}
and assume that it satisfies assumptions \eqref{eq:ass_drift_super} and \eqref{eq:ass_phi}.
Then, if the ODE \eqref{eq:ODE} with $K \equiv 0$ has an increasing solution $h_\nu(x_n)$ then 
$$
u(x) = h_\nu(x_n)
$$
solves the corresponding equation \eqref{eq:main}. 
To see this, note that its Hessian has one eigenvalue
$h_\nu''$ and $n-1$ zero eigenvalues. 
If $\Phi\geq0$, then $h_\nu''=-\Phi/\lambda\leq0$, 
while if $\Phi\leq0$, then $h_\nu''=-\Phi/\Lambda\geq0$. 
In either case, the equation follows directly from the ODE \eqref{eq:ODE} with $K\equiv0$.

Since $h_\nu$ is increasing, this example satisfies
$$
M(R) = h_\nu(R), \quad \text{and} \quad M'(R)=f_\nu(R).
$$
Thus, whenever the right-hand side of the estimate in
Theorem \ref{th:Phragmen--Lindelof} is positive, 
the resulting lower growth
estimate is optimal up to a multiplicative constant.
If 
$$
\frac{f_{\nu,R}(R)}{f_\nu(R)} \rightarrow 1,
$$ 
then the derivative estimate is also sharp in its multiplicative constant.

Theorem \ref{th:Phragmen--Lindelof} improves \cite[Theorem 2.1]{L22} in the following three ways: 
$(i)$ The function $\Phi$ controlling the lower-order terms is allowed to depend on the solution itself, 
$(ii)$ We allow for blow-up solutions, that is, solutions that approach infinity at a finite distance ($R_\infty < \infty$) from the boundary, 
$(iii)$ Rather than verifying whether a given PDE satisfies a set of technical assumptions, it is enough in the present version of the theorem to find growing positive solutions of the ODEs \eqref{eq:ODE}. 
In Section \ref{sec:applic} we consider several explicit examples in cases where \cite{L22} does not apply.

\bigskip

Related work concerns Phragm\'en--Lindel\"of principles,
Liouville theorems, and the influence of lower-order terms
on global solvability and blow-up. In the half-space setting,
Braga \cite{B17} and Lian \cite{L26} establish rigidity and
Liouville results for fully nonlinear elliptic equations.
Fuentes and Quaas \cite{FQ22} prove one-dimensional symmetry
for Hamilton--Jacobi equations involving Pucci operators,
while Filippucci, Pucci and Souplet \cite{FPS20} connect
half-space classification results with sharp gradient
blow-up profiles.

A closely related approach based on radial ODE analysis
is developed by Capuzzo Dolcetta, Leoni and Vitolo
\cite{CLV16}, who characterize the existence of entire
subsolutions and derive a priori bounds under generalized
Keller--Osserman conditions. Alarc\'on, Garc\'ia-Meli\'an
and Quaas \cite{AGQ12} study boundary blow-up problems
involving both reaction and gradient terms.
Further results include Liouville theorems for partial
trace equations with nonlinear gradient terms
by  Kindu, Mohammed and Tsegaw \cite{KMT25}, 
generalized Keller--Osserman criteria for
entire subsolutions of $k$-Hessian type equations with
gradient terms by Ji, Jiang and Li \cite{JJL23}, 
and local boundary lower estimates obtained through 
barrier constructions for convex solutions in bounded 
convex domains by Li and Li \cite{LL23}. 
Dancer, Daners and Hauer \cite{DDH15} establish Liouville
theorems for $p$-harmonic functions in exterior domains
and classify the behaviour at infinity of semibounded
solutions with homogeneous Neumann boundary conditions.
Their results concern asymptotic classification in exterior
domains, whereas we obtain lower growth estimates for
subsolutions in half-spaces.

Besides Lundstr\"om \cite{L22}, the reduction of PDE estimates to one-dimensional
differential inequalities also underlies two related
works. Fabricius and Lundstr\"om \cite{FL26} establish
sharp pointwise upper decay estimates for viscosity
subsolutions in thin domains, making explicit the
interaction between the transverse geometry and the
ellipticity, which may degenerate at infinity.
Lundstr\"om, Olofsson and Singh \cite{LOS25} develop
boundary estimates allowing ellipticity to degenerate
near the boundary and permitting general lower-order
terms. Applications include boundary Harnack
inequalities and H\"older continuity of quotients of positive solutions vanishing near $C^{1,1}$-boundaries.

Phragm\'en--Lindel\"of principles have also been established in other geometries and for different boundary conditions. 
Related results for the
$p$-Laplace equation in sectors and domains associated
with lower-dimensional hyperplanes can be found in
\cite{L16,LS23,LOS25}.
Abiko \cite{A23} treats fully nonlinear elliptic equations
on epigraphs and slabs with nonlinear dynamical boundary
conditions involving spatial derivatives and a time
derivative. Lamberti and Moroz \cite{LM24} establish a
superposition principle and a Phragm\'en--Lindel\"of
comparison principle for $p$-Laplace type equations in
radially symmetric domains, recovering and improving
lower and upper estimates for sub- and supersolutions
in the presence of Hardy-type potentials.
Fletcher and Steranka \cite{FS25} allow logarithmic
boundary growth in a Phragm\'en--Lindel\"of principle
for sub-$F$-extremals in sectors and apply it to Julia
limiting directions of quasiregular mappings.
Juraeva \cite{J24} obtains corresponding results for
biharmonic functions using Carleman-type integral formulas.
Broader applications include convergence results for
elliptic Dirichlet problems with variable coefficients
 by Peng \cite{P24a}, existence and multiplicity results for
singular fractional-order problems by Peng \cite{P24b}, and
small-data global existence and decay estimates for
fluid systems in infinite cylinders by Wu \cite{W23}.

The present paper focuses on lower growth bounds for
subsolutions in half-spaces, characterized through
associated ODEs. In particular, the method accommodates
solution- and gradient-dependent terms and describes
both growth at infinity and blow-up at a finite value
of the domain parameter.

%%%%%%%%%%%%%%%%%%%%%%%%%%%%%%%%%%%%%%%%%%%%%%%%%%%%%%%%%%%%%%%%%%%%%
%%%%%%%%%%%%%%%%%%%%%%%%%%%%%%%%%%%%%%%%%%%%%%%%%%%%%%%%%%%%%%%%%%%%%
%%%%%%%%%%%%%%%%%%%%%%%%%%%%%%%%%%%%%%%%%%%%%%%%%%%%%%%%%%%%%%%%%%%%%
%%%%%%%%%%%%%%%%%%%%%%%%%%%%%%%%%%%%%%%%%%%%%%%%%%%%%%%%%%%%%%%%%%%%%
%%%%%%%%%%%%%%%%%%%%%%%%%%%%%%%%%%%%%%%%%%%%%%%%%%%%%%%%%%%%%%%%%%%%%
%%%%%%%%%%%%%%%%%%%%%%%%%%%%%%%%%%%%%%%%%%%%%%%%%%%%%%%%%%%%%%%%%%%%%
%%%%%%%%%%%%%%%%%%%%%%%%%%%%%%%%%%%%%%%%%%%%%%%%%%%%%%%%%%%%%%%%%%%%%
%%%%%%%%%%%%%%%%%%%%%%%%%%%%%%%%%%%%%%%%%%%%%%%%%%%%%%%%%%%%%%%%%%%%%
%%%%%%%%%%%%%%%%%%%%%%%%%%%%%%%%%%%%%%%%%%%%%%%%%%%%%%%%%%%%%%%%%%%%%
%%%%%%%%%%%%%%%%%%%%%%%%%%%%%%%%%%%%%%%%%%%%%%%%%%%%%%%%%%%%%%%%%%%%%
%%%%%%%%%%%%%%%%%%%%%%%%%%%%%%%%%%%%%%%%%%%%%%%%%%%%%%%%%%%%%%%%%%%%%
%%%%%%%%%%%%%%%%%%%%%%%%%%%%%%%%%%%%%%%%%%%%%%%%%%%%%%%%%%%%%%%%%%%%%
%%%%%%%%%%%%%%%%%%%%%%%%%%%%%%%%%%%%%%%%%%%%%%%%%%%%%%%%%%%%%%%%%%%%%
%%%%%%%%%%%%%%%%%%%%%%%%%%%%%%%%%%%%%%%%%%%%%%%%%%%%%%%%%%%%%%%%%%%%%
%%%%%%%%%%%%%%%%%%%%%%%%%%%%%%%%%%%%%%%%%%%%%%%%%%%%%%%%%%%%%%%%%%%%%

\section{Proof of main results: Theorem \ref{th:Phragmen--Lindelof}}

\setcounter{theorem}{0}
\setcounter{equation}{0}

By $\Omega$ we denote a domain, that is, an open connected set.
For a set $E \subset \mathbb{R}^n$ we let $\overline{E}$ and $\partial E$ denote the closure and the boundary of $E$, respectively, and
points $x \in \mathbb{R}^n$ are sometimes denoted as $x = (x',x_n)$. 
By $c$ we denote a positive constant, not necessarily the same at each occurrence, 
we write $A \precsim B$ if there exists $c$ such that $A \leq cB$, and  $A(R)\asymp B(R)$ means that there exists $c>0$ such that
$c^{-1}B(R)\leq A(R)\leq cB(R)$ when $R\to R_\infty$.

A function $u : \Omega \to \mathbb{R}$ is a classical subsolution (supersolution) to \eqref{eq:main} in $\Omega$ if it is twice differentiable in $\Omega$ and satisfies $F(x,u,Du,D^2u) \leq 0$ ($\geq 0$). If the inequality is strict, $u$ is said to be a \textit{strict} classical subsolution (supersolution). A solution is a function that is both a sub- and a supersolution.

We choose to present our main results for viscosity solutions as defined in Crandall, Ishii and Lions \cite{CIL92} for $F:\mathbb{R}^n \times \mathbb{R} \times \mathbb{R}^n  \times \mathbb{S}^n \rightarrow \mathbb{R}$ a continuous function. However, neither the choice of viscosity solutions nor the assumption that $F$ is continuous are necessary for our results. In fact, other definitions of ``weak solutions" to \eqref{eq:main} can be considered whenever more appropriate for the equation at hand, as long as such weak subsolutions are upper semi-continuous (USC) and can be compared to classical strict supersolutions. In particular, we do not need a general comparison principle for sub- and supersolutions of \eqref{eq:main}, but only the following simple lemma which is immediate in the context of viscosity solutions:
\begin{lemma}\label{le:comp-weak-some-sence}
Suppose that the monotonicity condition in \eqref{eq:ass_drift_super} holds.
Let $\Omega$ be a bounded domain, $u \in USC(\overline{\Omega})$ a weak subsolution and $v$ a classical strict supersolution to \eqref{eq:main} in $\Omega$, continuous on $\overline{\Omega}$. If $u \leq v$ on $\partial \Omega$ then $u \leq v$ in $\Omega$.
\end{lemma}

The proof of Theorem \ref{th:Phragmen--Lindelof} relies on the construction of a classical strict supersolution of \eqref{eq:main} and comparison using Lemma \ref{le:comp-weak-some-sence}. The proof is in many parts similar to the proof of Lundstr\"om \cite[Theorem 2.1]{L22}. 
The following Lemma gives the desired strict classical supersolution:

%%%%%%%%%%%%%%%%%%%%%%%%%%%%%%%%%%%%%%%%%%%%
%%%%%%%%%%%%%%%%%%%%%%%%%%%%%%%%%%%%%%%%%%%%
%%%%%%%%%%%%%%%%%%%%%%%%%%%%%%%%%%%%%%%%%%%%

\begin{lemma}
\label{le:barrier_super}
Suppose that \eqref{eq:ass_drift_super} and \eqref{eq:ass_phi} hold, 
let $R > 0$ and define
$$
\Xi_R(x) = %\sqrt{ \sum_{i=1}^{n-1} x_i^2 + (x_n + \gamma)^2} - \gamma = 
|(x', x_n+\gamma)| - \gamma
$$
in which $\gamma = \gamma(R)$ is as in \eqref{eq:domain_2}. If the initial value problem \eqref{eq:ODE} has a solution $h_{\nu,R}' > 0$ on $(0,R) $, then the function
\begin{align*}%\label{eq:super}
V(x) =   h_{\nu,R}(\Xi_R(x))
\end{align*}
is a strict classical supersolution of \eqref{eq:main} in $D(R)$.
\end{lemma}

%%%%%%%%%%%%%%%%%%%%%%%%%%%%%%%%%%%%%%%%%%%%
%%%%%%%%%%%%%%%%%%%%%%%%%%%%%%%%%%%%%%%%%%%%
%%%%%%%%%%%%%%%%%%%%%%%%%%%%%%%%%%%%%%%%%%%%

\noindent
{\bf Proof.} %of Lemma \ref{le:barrier_super}.}
First note that, by construction, we have $|x| \geq \Xi_R(x) \geq x_n$. 
For notational simplicity we set
$\Xi = \Xi_R(x)$ and $h(t) = h_{\nu,R}(t)$.  
Differentiating yields
\begin{align*}
\frac{\partial V}{\partial x_i} = \frac{x_i}{|(x', x_n+\gamma)|} h'\left(\Xi\right), \quad \mbox{for} \quad 1 \leq i \leq n-1, \qquad
\frac{\partial V}{\partial x_n} = \frac{x_n + \gamma}{|(x', x_n+\gamma)|} h'\left(\Xi\right)
\end{align*}
and thus
$
\vert D V \vert = h'\left(\Xi\right).
$
The second derivatives become
\begin{align*}
\frac{\partial^2 V}{\partial x_i^2}
= %&
 \left( \frac{x_i}{|( x', x_n+\gamma)|} \right)^2 h''(\Xi)%\\
+%&
  \left(\frac{1}{|(x', x_n+\gamma)|} - \frac{x_i^2}{|(x', x_n+\gamma)|^3} \right) h'\left(\Xi\right),
  % \quad 1 < i < n, \quad\\
%
\end{align*}
for $1 \leq i \leq n-1$, and
\begin{align*}
\frac{\partial^2 V}{\partial x_n^2} &=  \left(\frac{x_n + \gamma}{|(x', x_n+\gamma)|}\right)^2 h''(\Xi)% \\
+%&
 \left(\frac{1}{|(x', x_n+\gamma)|} - \frac{(x_n + \gamma)^2}{|(x', x_n+\gamma)|^3} \right) h'\left(\Xi\right),
\end{align*}
giving
\begin{align*}
\text{Tr}(D^2 V) %&
=  h''(\Xi) + \frac{n-1}{|(x', x_n+\gamma)|}h'\left(\Xi\right).
\end{align*}

Assume now that $\Phi\geq 0$. 
By construction we then have from \eqref{eq:ODE} that 
$$h''(t) =  -  \frac{\Phi(t,h(t),h'(t))}{\lambda(t)} - K \frac{\Lambda(t)}{\lambda(t)} h'(t)
$$
and hence, since $V$ is radial, 
we can decompose $D^2 V = \left(D^2 V \right)^+ - \left(D^2 V \right)^-$ so that 
\begin{align*}
\text{Tr} \left(D^2 V \right)^+ &=  \frac{n-1}{|(x', x_n+\gamma)|}h'\left(\Xi\right) \quad \textrm{and} \\%\quad
\text{Tr} \left(D^2 V \right)^- &=  \frac{\Phi\left(\Xi, h\left(\Xi\right), h'\left(\Xi\right)\right)}{\lambda(\Xi)} + K \frac{\Lambda(\Xi)}{\lambda(\Xi)} h'(\Xi).
\end{align*}
Utilizing  the structure assumption \eqref{eq:ass_drift_super}, the fact that $V = h\left(\Xi\right)$ and $\vert D V \vert = h'\left(\Xi\right)$
we get
\begin{align}\label{eq:super_proof_2}
F(x,V,DV,D^2V) \geq & - \Phi\left(t, h\left(\Xi\right), h'\left(\Xi\right)\right) - \Lambda(x_n) \frac{n-1}{|(x', x_n+\gamma)|}h'\left(\Xi\right)\notag\\
&+ \frac{\lambda(x_n)}{\lambda(\Xi)}\Phi\left(\Xi,h\left(\Xi\right), h'\left(\Xi\right)\right) + K \frac{\lambda(x_n)\Lambda(\Xi)}{\lambda(\Xi)} h'(\Xi)   \notag\\
\geq  &- \Lambda(x_n)  \frac{n-1}{|(x', x_n+\gamma)|}h'(\Xi)
+   K \frac{\lambda(x_n)\Lambda(\Xi)}{\lambda(\Xi)} h'(\Xi).
\end{align}
For the last inequality, we used 
\begin{align}\label{eq:clarify_monotonicity_assumptions}
\frac{\lambda(x_n)}{\lambda(\Xi)} \geq 1, \quad \Phi \geq 0, \quad \text{and} \quad \Phi(t, \cdot, \cdot) \leq \Phi(\Xi, \cdot, \cdot),
\end{align}
which follows by the assumptions on $\lambda,  \Lambda$ and $\Phi$ in \eqref{eq:ass_drift_super} and \eqref{eq:ass_phi}, 
and by the geometric fact that
$$
x_n \leq \Xi(x) \leq |x|.
$$
In particular, 
the first inequality in \eqref{eq:clarify_monotonicity_assumptions} follows 
since $x_n \leq \Xi(x)$ and $\lambda$ is nonincreasing.
The last inequality follows since if $t = x_n$, then $\Phi$ is  nondecreasing, while if $t = |x|$, then $\Phi$ is nonincreasing, recall assumption \eqref{eq:ass_phi}.

To prove the supersolution property we need $F(x,V,DV,D^2V) > 0$ and since $h'(\Xi) > 0$ by assumption it thus suffices to ensure
\begin{align}\label{eq:lala}
\frac{n-1}{|(x', x_n+\gamma)|} <  K \frac{\lambda(x_n)\Lambda(\Xi)}{\Lambda(x_n)\lambda(\Xi)}.
\end{align}
Since $\lambda/\Lambda$ is nonincreasing, we have $\frac{\lambda(x_n)}{\Lambda(x_n)} \geq \frac{\lambda(\Xi)}{\Lambda(\Xi)}$ and it suffices that
$$
\frac{n-1}{|(x', x_n+\gamma)|} <  K .
$$
However, as $\gamma \leq |(x', x_n+\gamma)|$ in $\mathbb{R}_+^n$ it is even enough to have
$$
\frac{n-1}{\gamma(R)} <  K
$$
which holds by construction since $\gamma >0$ and $K = {n}/{\gamma}$.

If $\Phi \leq 0$, then
$$
h''(t) = -
\frac{\Phi\left(t,h\left(t\right),h'\left(t\right)\right)}{\Lambda(t)} -  K \frac{\Lambda(t)}{\lambda(t)} h'(t)
$$
and we obtain, after decomposing the Hessian, that
\begin{align*}
\operatorname{Tr}(D^2V)^+=(A-B)^+
+\frac{n-1}{|(x', x_n+\gamma)|}h'(\Xi),
\quad
\operatorname{Tr}(D^2V)^-=(A-B)^-,
\end{align*}
in which $t^\pm=\max\{\pm t,0\}$, 
$$
A=-\frac{\Phi(\Xi,h(\Xi),h'(\Xi))}{\Lambda(\Xi)}\geq 0, \quad \text{and}\quad B=K\frac{\Lambda(\Xi)}{\lambda(\Xi)}h'(\Xi) > 0.
$$
Using that $\Lambda (A-B)^+-\lambda (A-B)^-\leq \Lambda A-\lambda B$ we end up with
\begin{align}\label{eq:super_proof_33}
F(x,V,DV,D^2V)
\geq&- \Phi\left(t, h\left(\Xi\right), h'\left(\Xi\right)\right) -\Lambda(x_n) \frac{n-1}{|(x', x_n+\gamma)|} h'(\Xi) \notag\\
&- \Lambda(x_n) \left(A - B\right)^+  + \lambda(x_n)  \left(A - B\right)^-\notag\\
\geq&- \Phi\left(t, h\left(\Xi\right), h'\left(\Xi\right)\right) -\Lambda(x_n) \frac{n-1}{|(x', x_n+\gamma)|} h'(\Xi) \notag\\
&- \Lambda(x_n) A  + \lambda(x_n)  B\notag\\
=&- \Phi\left(t, h\left(\Xi\right), h'\left(\Xi\right)\right) -\Lambda(x_n) \frac{n-1}{|(x', x_n+\gamma)|} h'(\Xi) \notag\\
& + \frac{\Lambda(x_n)}{\Lambda(\Xi)} \Phi\left(\Xi,h\left(\Xi\right),h'\left(\Xi\right)\right)  + \lambda(x_n) K \frac{\Lambda(\Xi)}{\lambda(\Xi)} h'(\Xi)\notag\\
\geq &-\Lambda(x_n)\frac{n-1}{|(x', x_n+\gamma)|} h'(\Xi) + \lambda(x_n) K \frac{\Lambda(\Xi)}{\lambda(\Xi)} h'(\Xi).
\end{align}
For the last inequality, we used 
\begin{align*}
\frac{\Lambda(x_n)}{\Lambda(\Xi)} \leq 1, \quad \Phi \leq 0, \quad \text{and} \quad -\Phi(t, \cdot, \cdot) \geq -\Phi(\Xi, \cdot, \cdot),
\end{align*}
which again follows by the assumptions on $\lambda,  \Lambda, \Phi$ in \eqref{eq:ass_drift_super} and \eqref{eq:ass_phi}, 
and by $x_n \leq \Xi(x) \leq |x|.$ 
Indeed, the first inequality follows since $x_n \leq \Xi(x)$ and $\Lambda$ is nondecreasing, 
while the last inequality is the same as the last inequality in \eqref{eq:clarify_monotonicity_assumptions}.

As $h'>0$ by assumption, what remains to ensure that $V$ is a strict supersolution is to show that
$$
\frac{n-1}{|(x', x_n+\gamma)|} < K \frac{\lambda(x_n)}{\Lambda(x_n)} \frac{\Lambda(\Xi)}{\lambda(\Xi)}.
$$
We are thus back in the same situation as in \eqref{eq:lala} and the proof of Lemma \ref{le:barrier_super} is complete.
$\hfill \Box$ \\

%%%%%%%%%%%%%%%%%%%%%%%%%%%%%%%%%%%%%%%%%%%%
%%%%%%%%%%%%%%%%%%%%%%%%%%%%%%%%%%%%%%%%%%%%
%%%%%%%%%%%%%%%%%%%%%%%%%%%%%%%%%%%%%%%%%%%%

\noindent
{\bf Proof of Theorem \ref{th:Phragmen--Lindelof}.}
Let $\nu > 0$ be as in the theorem and apply Lemma \ref{le:barrier_super} to ensure, 
for any $R \in (0, R_\infty)$, 
that $V := h_{\nu,R}(\Xi_R(x))$ is a strict supersolution in $D(R)$.
We define
$$
C=M(R)-h_{\nu,R}(R),
$$
and prove that then we must have $C \geq 0$. 
Assume $C < 0$ and observe that then
$$
u \leq M(R) < V \quad \text{on} \quad \partial D(R) \cap \mathbb{R}^n_+.
$$
Since $V \geq 0$ on $\partial D(R) \cap \partial \mathbb{R}^n_+$  
the boundary condition on $u$ yields 
\begin{align*}
\limsup_{x \to z} u(x) \leq V(z) \quad \textrm{for all} \quad z \in \partial D(R),
\end{align*}
and so by the weak comparison principle in Lemma \ref{le:comp-weak-some-sence} we get $u \leq V$ in $D(R)$.
By the assumptions on the given subsolution $u$, $h_\nu$ and $h_{\nu,R}$ it follows that, for $\bar x_n < R < R_\infty$, 
$$
h_{\nu}(\bar x_n)  \leq  u(\bar x) \leq V(\bar x) = h_{\nu,R}(\bar x_n)
<h_{\nu}(\bar x_n),
$$
which gives a contradiction. 
Thus, we conclude that $C\geq 0$, and by the monotonicity assumption in \eqref{eq:ass_drift_super} 
it follows that $V + C$ remains a strict supersolution of  \eqref{eq:main}.
Moreover, we have $V + C = M(R)$ on $\partial D(R) \cap \mathbb{R}^n_+$ and thus, 
by the comparison principle in Lemma \ref{le:comp-weak-some-sence} we conclude that $V + C \geq u$ in $D(R)$.

Now, for $\bar{x}_n<R<R_\infty$ and $\delta \in (0,R)$ it follows, since $\gamma(R)$ is nondecreasing and $f_{\nu,R} = h_{\nu,R}'>0$, that
\begin{align*}
M(R) - M(R- \delta) &\geq V(0, \dots, 0, R) - V(0, \dots, 0, R - \delta)\\
&= h_{\nu,R}(R) - h_{\nu,R}(R-\delta) > 0, 
\end{align*}
for every $ \delta \in (0,R)$. 
Thus, $M(R)$ is increasing on $(\bar x_n, R_\infty)$ 
and by taking the limit we see that
\begin{align*}
M'(R) = \liminf_{\delta \to 0^+} \frac{M(R) - M(R-\delta)}{ \delta}
\geq
\liminf_{\delta \to 0^+} \frac{h_{\nu,R}(R) - h_{\nu,R}(R-\delta) }{ \delta} = f_{\nu,R}(R).
\end{align*}
Hence
$$
\liminf_{R \to R_\infty}\frac{M'(R)}{f_{\nu}\left(R\right)} \geq \liminf_{R \to R_\infty}\frac{f_{\nu,R}\left(R\right)}{f_{\nu}\left(R\right)}.
$$
This completes the proof of the theorem. $\hfill \Box$\\

%%%%%%%%%%%%%%%%%%%%%%%%%%%%%%%%%%%%%%%%%%%%%%%%%%%%%%%%%%%%%%%%%%%%%
%%%%%%%%%%%%%%%%%%%%%%%%%%%%%%%%%%%%%%%%%%%%%%%%%%%%%%%%%%%%%%%%%%%%%
%%%%%%%%%%%%%%%%%%%%%%%%%%%%%%%%%%%%%%%%%%%%%%%%%%%%%%%%%%%%%%%%%%%%%
%%%%%%%%%%%%%%%%%%%%%%%%%%%%%%%%%%%%%%%%%%%%%%%%%%%%%%%%%%%%%%%%%%%%%
%%%%%%%%%%%%%%%%%%%%%%%%%%%%%%%%%%%%%%%%%%%%%%%%%%%%%%%%%%%%%%%%%%%%%
%%%%%%%%%%%%%%%%%%%%%%%%%%%%%%%%%%%%%%%%%%%%%%%%%%%%%%%%%%%%%%%%%%%%%
%%%%%%%%%%%%%%%%%%%%%%%%%%%%%%%%%%%%%%%%%%%%%%%%%%%%%%%%%%%%%%%%%%%%%
%%%%%%%%%%%%%%%%%%%%%%%%%%%%%%%%%%%%%%%%%%%%%%%%%%%%%%%%%%%%%%%%%%%%%
%%%%%%%%%%%%%%%%%%%%%%%%%%%%%%%%%%%%%%%%%%%%%%%%%%%%%%%%%%%%%%%%%%%%%
%%%%%%%%%%%%%%%%%%%%%%%%%%%%%%%%%%%%%%%%%%%%%%%%%%%%%%%%%%%%%%%%%%%%%
%%%%%%%%%%%%%%%%%%%%%%%%%%%%%%%%%%%%%%%%%%%%%%%%%%%%%%%%%%%%%%%%%%%%%
%%%%%%%%%%%%%%%%%%%%%%%%%%%%%%%%%%%%%%%%%%%%%%%%%%%%%%%%%%%%%%%%%%%%%
%%%%%%%%%%%%%%%%%%%%%%%%%%%%%%%%%%%%%%%%%%%%%%%%%%%%%%%%%%%%%%%%%%%%%
%%%%%%%%%%%%%%%%%%%%%%%%%%%%%%%%%%%%%%%%%%%%%%%%%%%%%%%%%%%%%%%%%%%%%
%%%%%%%%%%%%%%%%%%%%%%%%%%%%%%%%%%%%%%%%%%%%%%%%%%%%%%%%%%%%%%%%%%%%%

\section{Applications of Theorem \ref{th:Phragmen--Lindelof}}
\label{sec:applic}

\setcounter{theorem}{0}
\setcounter{equation}{0}

In this section we apply Theorem \ref{th:Phragmen--Lindelof} in some cases where we can find or estimate solutions of the ODEs in \eqref{eq:ODE}, calculate the limit
\begin{align}\label{eq:limit_appl}
\liminf_{R\to R_\infty}\frac{f_{\nu,R}(R)}{f_{\nu}(R)},
\end{align}
and thereby conclude explicit growth estimates. 
Our examples show that the assumptions are fairly permitting, allowing for several interesting PDEs not covered by Lundstr\"om in \cite{L22}.

%%%%%%%%%%%%%%%%%%%%%%%%%%%%%%%%%%%%%%%%%%%%%%%%%%%%%%%%%%%%%%%%%%%%%
%%%%%%%%%%%%%%%%%%%%%%%%%%%%%%%%%%%%%%%%%%%%%%%%%%%%%%%%%%%%%%%%%%%%%
%%%%%%%%%%%%%%%%%%%%%%%%%%%%%%%%%%%%%%%%%%%%%%%%%%%%%%%%%%%%%%%%%%%%%
%%%%%%%%%%%%%%%%%%%%%%%%%%%%%%%%%%%%%%%%%%%%%%%%%%%%%%%%%%%%%%%%%%%%%
%%%%%%%%%%%%%%%%%%%%%%%%%%%%%%%%%%%%%%%%%%%%%%%%%%%%%%%%%%%%%%%%%%%%%

\subsection{The case $\Phi(t,u,|Du|) =  - C(t) \leq 0$}
In this case 
the initial value problem \eqref{eq:ODE} simplifies to 
$$
f'(t) = \frac{C(t)}{\Lambda(t)} - K(R) \frac{\Lambda(t)}{\lambda(t)} f(t),   \quad t \in (0,R) \quad \text{with} \quad h(0)= 0,\quad f(0) = \nu > 0.
$$ 
Calculations give
\begin{align}\label{eq:case1-f1}
f_{\nu,R} (t) = e^{-K(R) \int_{0}^{t} \frac{\Lambda(s)}{\lambda(s)}ds} \left( \nu + \int_{0}^{t} e^{K(R) \int_{0}^{\bar t} \frac{\Lambda(s)}{\lambda(s)}ds } \frac{C(\bar t)}{\Lambda(\bar t)} d\bar t \right),
\end{align}
and
\begin{align}\label{eq:case1-f2}
f_{\nu} (t) =  \nu + \int_{0}^{t}  \frac{C(\bar t)}{\Lambda(\bar t)} d\bar t.
\end{align}
Note that these solutions exist for $t \in (0, R_\infty)$ with $R_\infty=\infty$.
The limit \eqref{eq:limit_appl} thus turns to
\begin{align*}
\liminf_{R\to\infty}\frac{e^{-K(R) \int_{0}^{R} \frac{\Lambda(s)}{\lambda(s)}ds} \left( \nu + \int_{0}^{R} e^{K(R) \int_{0}^{\bar t} \frac{\Lambda(s)}{\lambda(s)}ds } \frac{C(\bar t)}{\Lambda(\bar t)} d\bar t \right)}{\nu + \int_{0}^{R}  \frac{C(\bar t)}{\Lambda(\bar t)} d\bar t},
\end{align*}
which will be positive if 
\begin{align}\label{eq:gamma_first}
K(R) \int_{0}^{R} \frac{\Lambda(s)}{\lambda(s)} ds = \frac{n}{\gamma(R)} \int_{0}^{R} \frac{\Lambda(s)}{\lambda(s)} ds \precsim  1,
\end{align}
whenever $R \in \mathbb{R}_+$.
(As an example, in the case of uniform ellipticity $\Lambda / \lambda = constant$, this can be achieved by taking
$\gamma(R) = R$ so that $D(R)$ preserves the `standard geometry' of the classical
Phragmen–Lindelöf theorem.)

Whenever \eqref{eq:gamma_first} holds we obtain the following growth estimate for subsolutions of PDEs satisfying the assumptions of Theorem \ref{th:Phragmen--Lindelof}:
\begin{align}\label{eq:case1}
\liminf_{R\to\infty} \, \frac{M'(R)}{\nu + \int_{0}^{R}  \frac{C(\bar t)}{\Lambda(\bar t)} d\bar t} > 0.
\end{align}

\bigskip

\noindent
{\bf Example}. If $C(t)/\Lambda(t) = constant > 0$, then we retrieve the growth estimate
$$
\liminf_{R\to\infty} \, \frac{M'(R)}{R} > 0 \implies \liminf_{R\to\infty} \, \frac{M(R)}{R^2} > 0.
$$

\bigskip

\noindent
{\bf Example}. If $C(t)/\Lambda(t) = (t+1)^k$, $k\in\mathbb{R}$, 
then
$$
\int_{0}^{R}\frac{C(s)}{\Lambda(s)} ds = \left\{
\begin{array}{ll}
\frac{(R + 1)^{k + 1} - 1}{k + 1} & \text{if $k \neq -1$}, \\
\log(R + 1) & \text{if $k = -1$},
\end{array}
\right.
$$
and therefore
\begin{align*}
\liminf_{R\to\infty} \, M'(R) > 0 &\implies \liminf_{R\to\infty} \, \frac{M(R)}{R} > 0 \quad  \text{if $k < -1$}, \\
\liminf_{R\to\infty} \, \frac{M'(R)}{\log R} > 0 &\implies \liminf_{R\to\infty} \, \frac{M(R)}{R \log R} > 0 \quad  \text{if $k = -1$}, \\
\liminf_{R\to\infty} \, \frac{M'(R)}{R^{k + 1}} > 0 &\implies \liminf_{R\to\infty} \, \frac{M(R)}{R^{k + 2}} > 0 \quad  \text{if $k > -1$}.
\end{align*}

Concerning sharpness of \eqref{eq:case1} we note that if $\Phi(t,u,|Du|) =  - C(t) \leq 0$, then 
$$
u(x) = \int_{0}^{x_n} h'_\nu (t) dt = \nu\, x_n + \int_{0}^{x_n} \int_{0}^{s}  \frac{C(\bar t)}{\Lambda(\bar t)} d\bar t\, ds
$$
is a classical solution of \eqref{eq:main} with the operator \eqref{eq:sharp}, attaining the slowest growth allowed by \eqref{eq:case1}.

Two simple equations included in the setting of this subsection are
$\Delta u=|x|^k$ and $\Delta u=(1+x_n)^{-k}$, $k\geq0$,
corresponding to $\Phi=-|x|^k$, nonincreasing in $|x|$,
and $\Phi=-(1+x_n)^{-k}$, nondecreasing in $x_n$, respectively.
Finally, we note that all equations \eqref{eq:main} satisfying \eqref{eq:ass_drift_super} with $\Phi(t,u,|Du|) =  - C(t) \leq 0$ automatically satisfy \eqref{eq:ass_drift_super} with $\Phi \equiv 0$ and this case was outlined already in Lundstr\"om \cite{L22} but then only linear growth is obtained.

%%%%%%%%%%%%%%%%%%%%%%%%%%%%%%%%%%%%%%%%%%%%%%%%%%%%%%%%%%%%%%%%%%%%%
%%%%%%%%%%%%%%%%%%%%%%%%%%%%%%%%%%%%%%%%%%%%%%%%%%%%%%%%%%%%%%%%%%%%%
%%%%%%%%%%%%%%%%%%%%%%%%%%%%%%%%%%%%%%%%%%%%%%%%%%%%%%%%%%%%%%%%%%%%%
%%%%%%%%%%%%%%%%%%%%%%%%%%%%%%%%%%%%%%%%%%%%%%%%%%%%%%%%%%%%%%%%%%%%%
%%%%%%%%%%%%%%%%%%%%%%%%%%%%%%%%%%%%%%%%%%%%%%%%%%%%%%%%%%%%%%%%%%%%%

\subsection{The case $\Phi(t,u,|Du|) =  C(t) \geq 0$}
First, note that when $\Phi \geq 0$ the maximum principle does not hold in general as this source term 
may cause solutions to attain an interior maximum, making growth estimates void. 
Similar to the former case, the initial value problem \eqref{eq:ODE} simplifies to
$$
f'(t) = - \frac{C(t)}{\lambda(t)} - K(R) \frac{\Lambda(t)}{\lambda(t)} f(t),   \quad t \in (0,R) \quad \text{with} \quad h(0)= 0,\quad f(0) = \nu > 0,
$$
and therefore we have
\begin{align}\label{eq:case2-f1}
f_{\nu,R} (t) = e^{-K(R) \int_{0}^{t} \frac{\Lambda(s)}{\lambda(s)}ds} \left( \nu - \int_{0}^{t} e^{K(R) \int_{0}^{\bar t} \frac{\Lambda(s)}{\lambda(s)}ds } \frac{C(\bar t)}{\lambda(\bar t)} d\bar t \right),
\end{align}
\begin{align}\label{eq:case2-f2}
f_{\nu} (t) =  \nu - \int_{0}^{t}  \frac{C(\bar t)}{\lambda(\bar t)} d\bar t,
\end{align}
and the limit to investigate is thus
\begin{align*}
\liminf_{R\to\infty}\frac{f_{\nu,R}(R)}{f_{\nu}(R)} = \liminf_{R\to\infty}
\frac{e^{-K(R) \int_{0}^{R} \frac{\Lambda(s)}{\lambda(s)}ds} \left( \nu - \int_{0}^{R} e^{K(R) \int_{0}^{\bar t} \frac{\Lambda(s)}{\lambda(s)}ds } \frac{C(\bar t)}{\lambda(\bar t)} d\bar t \right)}{\nu - \int_{0}^{R}  \frac{C(\bar t)}{\lambda(\bar t)} d\bar t}.
\end{align*}
First, as in the former case, we demand that $\gamma$ satisfies \eqref{eq:gamma_first}, and secondly, we need
\begin{align}\label{eq:case2-solpos}
\int_{0}^{R}  \frac{C(s)}{\lambda(s)} ds < \nu \quad \text{whenever} \quad R \in \mathbb{R}_+
\end{align}
for $f_\nu$ to be positive.
%for all $t \in \mathbf{R}_+$.
As long as \eqref{eq:case2-solpos} holds, 
we can ensure positivity of $f_{\nu,R}$ by increasing $\gamma$ 
speeding up the convergence
$$
\frac{n}{\gamma(R)} \int_{0}^{R} \frac{\Lambda(s)}{\lambda(s)} ds \to 0  \quad \text{as} \quad R \to \infty
$$
and thereby forcing $e^{K(R) \int_{0}^{R} \frac{\Lambda(s)}{\lambda(s)}ds}$ to approach 1 fast enough as $R \to \infty$. 
Consequently, as long as \eqref{eq:case2-solpos} holds, we can produce a growth estimate by taking $\gamma$ large enough.

\bigskip

\begin{figure}
    \centering
    \includegraphics[height=6.0cm,width=7.5cm]{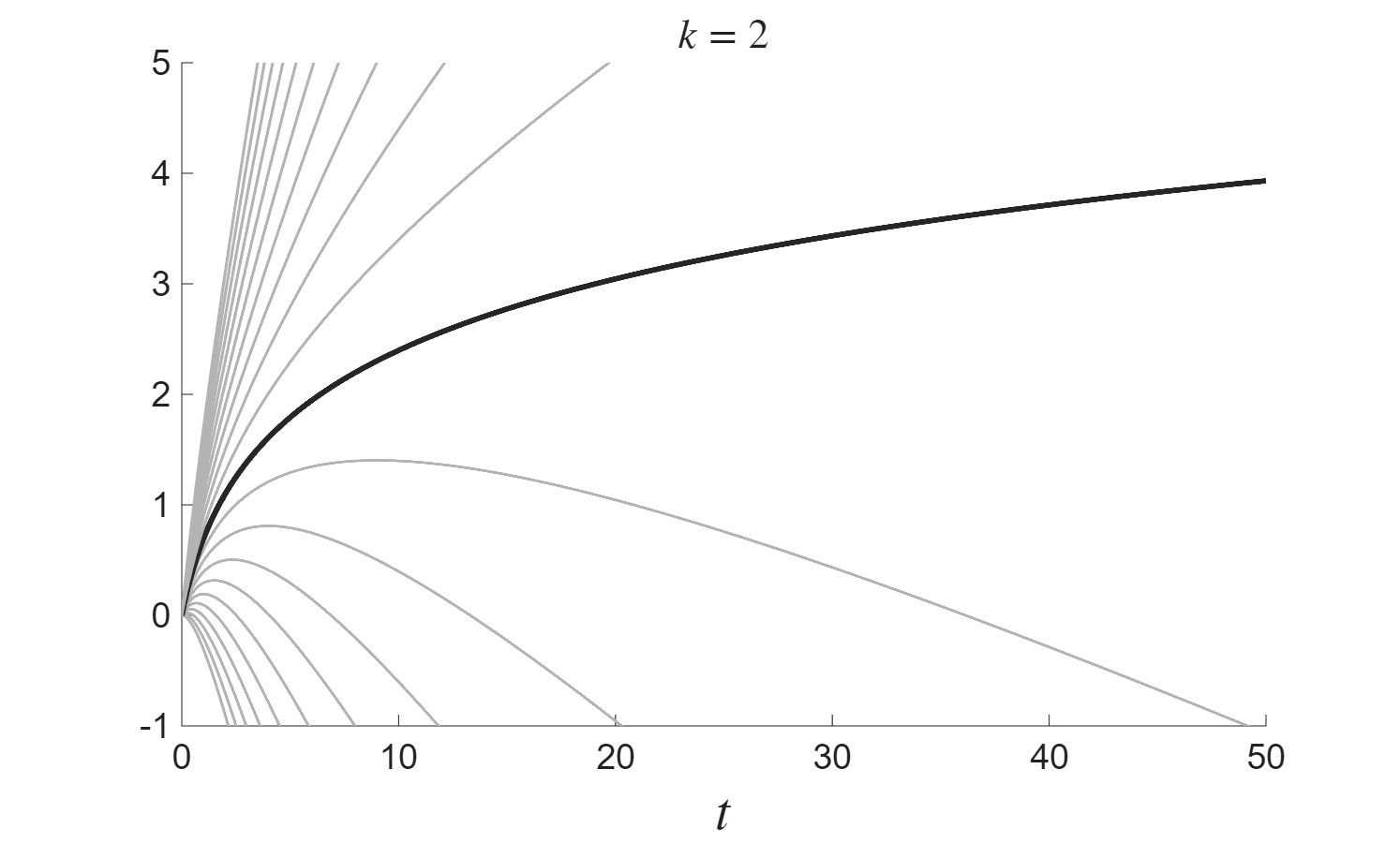}
    \includegraphics[height=6.0cm,width=7.5cm]{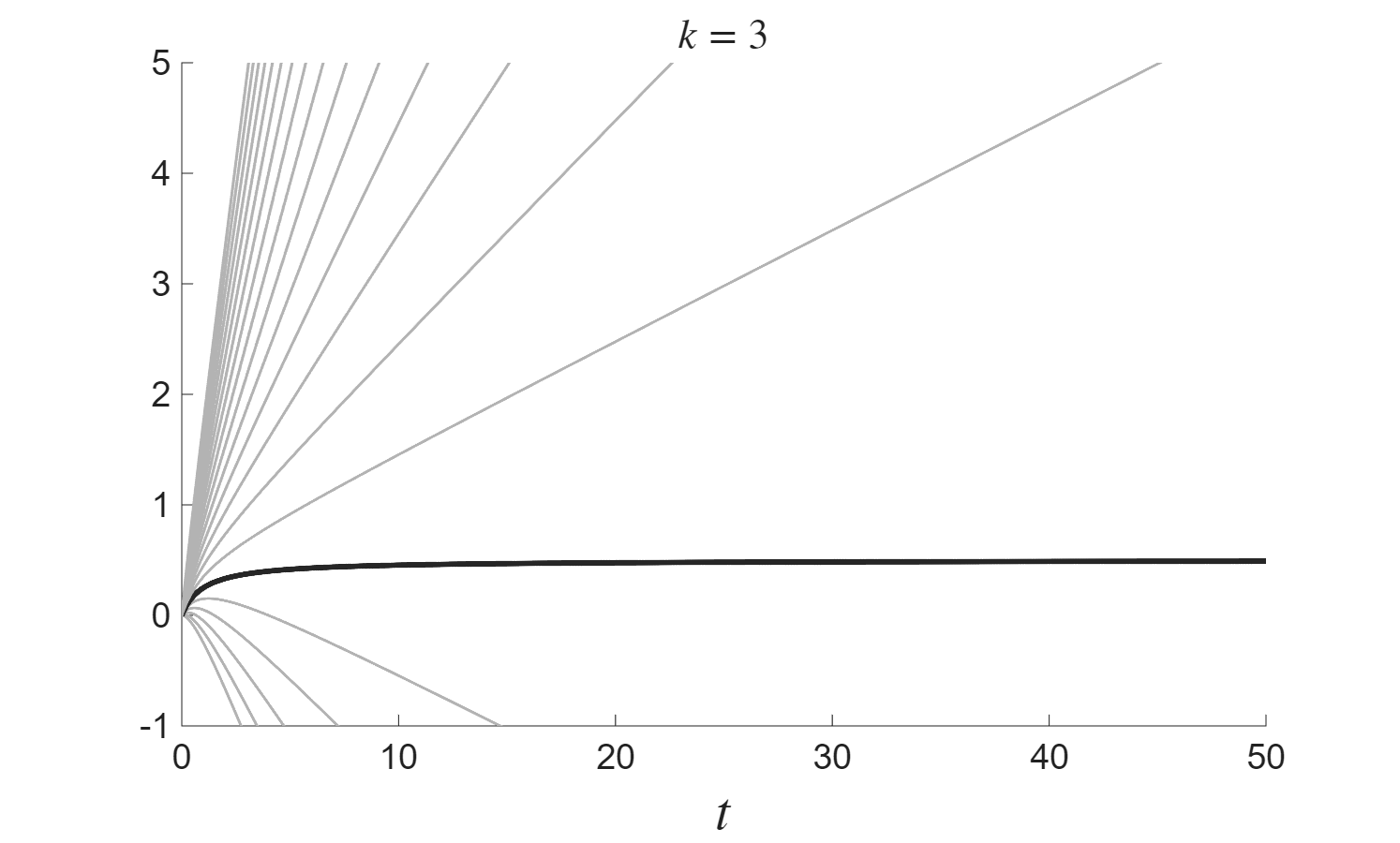}
    \caption{The function $h_\nu(t)$ for $\nu = [0,0.1,0.2, \dots, 2]$ and for $k=2$ (left) and $k=3$ (right). 
    The thick curves correspond to the borderline case $\nu = 1/(k-1)$.}
    \label{fig:1}
\end{figure}

\noindent
{\bf Example}. If $C(s)/\lambda(s) = \frac{1}{(1 + s)^k}$ with $k > 1$,
then we have
$$
f_\nu (R) = \nu - \int_0^R \frac{1}{(1+s)^k}\,ds \searrow\nu - \frac{1}{k-1} \quad \text{as} \quad R \to \infty,
$$
and conclude the following from Theorem \ref{th:Phragmen--Lindelof}:
\begin{align*}
\nu < \frac{1}{k-1}& \quad \Longrightarrow \quad \text{No growth estimate,}\\
\nu = \frac{1}{k-1}& \quad \Longrightarrow \quad \liminf_{R\to\infty} \, \frac{M'(R)}{R^{1-k}} > 0,\\
\nu > \frac{1}{k-1}& \quad \Longrightarrow \quad \liminf_{R\to\infty} \, \frac{M(R)}{R} > 0. 
\end{align*}
Thus, the slowest growing subsolution depends heavily on $\nu$, i.e. the value of the subsolution a bit inside the domain.
In the borderline case $\nu = \frac{1}{k-1}$ subsolutions may be bounded or must grow to infinity, depending on the value of $k$.
Indeed, from \eqref{eq:hejdutydligagrowthest} and
\begin{equation*}
h_\nu(R)=\int_0^R f_\nu(t)\,dt
=
\begin{cases}
\displaystyle
\left(\nu-\frac{1}{k-1}\right)R
+\frac{(1+R)^{2-k}-1}{(k-1)(2-k)},
& k\ne 2,\\%[6pt]
(\nu-1)R+\log(1+R),
& k=2,
\end{cases}
\end{equation*}
we realize that subsolutions may be bounded in the case $k > 2$. 
Otherwise, the borderline case implies
$$
\liminf_{R\to\infty} \, \frac{M(R)}{\log(R)} > 0 \quad \text{if $k = 2$} \quad \text{and} \quad
\liminf_{R\to\infty} \, \frac{M(R)}{R^{2-k}} > 0 \quad \text{if $1 < k < 2$}.
$$
Figure \ref{fig:1} shows $h_\nu$ for $k = 2$ and $k = 3$ and for several values of $\nu$ in each case. 
When $k = 2$ the borderline case happens at $\nu = 1$ giving logarithmic growth. All $\nu > 1$ give linear growth. 
When $k = 3$ the borderline case happens at $\nu = \frac{1}{2}$ allowing for bounded subsolutions. 
All $\nu > \frac{1}{2}$ give linear growth. 
A simple PDE included in this example is
$$
\Delta u + \frac{1}{(1 + |x|)^k} = 0
$$
and the reader can easily verify the above estimates for this PDE in the one-dimensional case.

%%%%%%%%%%%%%%%%%%%%%%%%%%%%%%%%%%%%%%%%%%%%%%%%%%%%%%%%%%%%%%%%%%%%%
%%%%%%%%%%%%%%%%%%%%%%%%%%%%%%%%%%%%%%%%%%%%%%%%%%%%%%%%%%%%%%%%%%%%%
%%%%%%%%%%%%%%%%%%%%%%%%%%%%%%%%%%%%%%%%%%%%%%%%%%%%%%%%%%%%%%%%%%%%%
%%%%%%%%%%%%%%%%%%%%%%%%%%%%%%%%%%%%%%%%%%%%%%%%%%%%%%%%%%%%%%%%%%%%%
%%%%%%%%%%%%%%%%%%%%%%%%%%%%%%%%%%%%%%%%%%%%%%%%%%%%%%%%%%%%%%%%%%%%%

\subsection{The case $\Phi(s,u,|Du|) =  -C(s)|Du|^k\leq 0$} \label{subsec:blowup}
We assume $k \in \mathbb{R}$ and that $C(t)$ is a continuous positive function for $t \in [0,\infty)$. 
This case will show blow-up solutions. 
For $k<0$, the notion of solution requires care because
of the singularity at $Du=0$. However, the proof only
requires the subsolution inequality for test functions
with nonzero gradient: the barriers constructed above
satisfy $|DV|=h'>0$, so the operator is evaluated only
at nonzero gradients in the comparison argument.
In this case the ODE in \eqref{eq:ODE} becomes
\begin{align*}%\label{eq:ODE}
\frac{df}{dt} =  \frac{C(t)}{\Lambda(t)}f^k(t) - K(R) \frac{\Lambda(t)}{\lambda(t)} f(t),  \quad t \in (0,R) \quad \text{with} \quad h(0)= 0,\quad f(0) = \nu > 0,
\end{align*}
%
%and we assume that $C(t)$ is a continuous positive function for $t \in [0,R]$. 
To simplify we replace the above equation with 
$$
\frac{d\tilde f}{dt} = A(t) \left( \tilde f^k - \widetilde{K} \tilde f \right), \quad t \in (0,R) \quad \text{with} \quad \tilde f(0) = \nu,
$$
where $A(t) = {C(t)}/{\Lambda(t)}$ and $\widetilde{K} = \widetilde{K}(R)$ is an upper bound of $\frac{K(R) \Lambda^2(t)}{C(t) \lambda(t)}$ for $t \in (0,R)$. This simplification is possible since solutions $\tilde f$ of this separable ODE will lie below solutions $f$ of the original ODE, thus producing a lower bound on the limit in Theorem \ref{th:Phragmen--Lindelof}. 
To find the solution for $k \neq 1$ we 
substitute \(z=\tilde f^{\,1-k}\) giving
\[
z'=(1-k)A(t)(1-\widetilde Kz),\quad z(0)=\nu^{1-k},
\]
which reveals the solution
\[
\tilde f_{\nu,R}(t)=
\begin{cases}
\displaystyle
\nu e^{(1-\widetilde K)\int_{0}^{t}A(s)ds },
& k=1,\\[6pt]
\displaystyle
\left[
\nu^{1-k}e^{(k-1)\widetilde K \int_{0}^{t}A(s)ds }
+\frac{1-e^{(k-1)\widetilde K \int_{0}^{t}A(s)ds }}{\widetilde K}
\right]^{\frac1{1-k}},
& k\ne1.
\end{cases}
\]
Similarly, solving $f'=A(t)f^k$ with $f(0)=\nu$ yields
\[
f_\nu(t)=
\begin{cases}
\displaystyle
\nu e^{\int_{0}^{t}A(s)ds },
& k=1,\\[6pt]
\displaystyle
\left[\nu^{1-k}+(1-k)\int_{0}^{t}A(s)ds \right]^{\frac1{1-k}},
& k\ne1.
\end{cases}
\]
For $k\leq1$, these solutions remain positive and finite
on every finite interval, and we set $R_\infty=\infty$.
For $k>1$, the solution $f_\nu$ explodes at a finite
$R_\infty$ if
\begin{align}\label{eq:R_infty}
 \int_{0}^{R_\infty}A(s)ds = \frac{\nu^{1-k}}{k - 1}.
\end{align}
If this level is not reached at a finite $R_\infty$,
then we set $R_\infty=\infty$.
Note that if $f_\nu$ explodes, it does so before $ \tilde f_{\nu,R}$.

Applying Theorem \ref{th:Phragmen--Lindelof} we thus conclude that, for $k = 1$,
\begin{align*}
\liminf_{R\to \infty} \, \frac{M'(R)}{h'_{\nu}(R)} \, \geq \,\liminf_{R\to \infty} \,   e^{ - \widetilde{K}\int_{0}^{R} A(s) ds},
\end{align*}
and for $k \neq 1$,
\begin{align}\label{eq:limit_hej}
\liminf_{R\to R_{\infty}} \, \frac{M'(R)}{h'_{\nu}(R)} \, \geq \,\liminf_{R\to R_{\infty}} \, 
\left(\frac{\displaystyle\nu^{1-k} + (1-k)I(R)}
{\displaystyle{\nu^{1-k}}e^{(k-1)\widetilde K I(R)}
+\frac{1 - e^{(k-1)\widetilde K I(R)}}{\widetilde K}} \right)^{\frac{1}{k-1}},
\end{align}
in which $I(R) = \int_{0}^{R} A(s)ds$. 
For $k = 1$ we conclude that choosing $\gamma = \gamma(R)$ such that
\begin{align}
    \widetilde{K}(R) \int_{0}^{R} A(s) ds \precsim 1 %\quad \text{as} \quad R \to \infty
\end{align}
ensures a growth estimate.
For \(k<1\), put \(a=(1-k)\widetilde K I(R)\) and observe that 
since \(1-e^{-a}\ge ae^{-a}\) we have
\[
\begin{aligned}
\tilde f_{\nu,R}(R)^{1-k}
&=\nu^{1-k}e^{-a}+\frac{1-e^{-a}}{\widetilde K}\\
&\ge e^{-a}\bigl[\nu^{1-k}+(1-k)I(R)\bigr]
=e^{-a}f_\nu(R)^{1-k}.
\end{aligned}
\]
Therefore, the same choice of $\gamma$ ensures a growth estimate also in this case.

For $k>1$, we rewrite the denominator in \eqref{eq:limit_hej} as
\[
e^{(k-1)\widetilde K I(R)}
\left[
\nu^{1-k}+(1-k)I(R)
+(k-1)\int_0^{I(R)}
\left(1-e^{-(k-1)\widetilde K s}\right)\,ds
\right].
\]
Since $1-e^{-x}\leq x$ for $x\geq0$, this is bounded above by
\[
e^{(k-1)\widetilde K I(R)}
\left[
\nu^{1-k}+(1-k)I(R)
+\frac{(k-1)^2\widetilde K I(R)^2}{2}
\right].
\]
Since $0\leq I(R)<\nu^{1-k}/(k-1)$, choosing
$\gamma=\gamma(R)$ such that
\begin{align}\label{eq:tildeKval}
\widetilde K(R)\precsim \nu^{1-k}+(1-k)I(R)
%\quad\text{as }R\to R_\infty
\end{align}
bounds the denominator by a constant multiple of the
numerator and hence ensures positivity of the lower
limit on the right-hand side of \eqref{eq:limit_hej}.

In case $C\equiv 1$ and $\lambda,\Lambda$ are constants, 
then
\[
I(R)=\int_0^R A(s)ds = \frac{R}{\Lambda}
\quad \text{and} \quad
R_\infty=\frac{\Lambda\nu^{1-k}}{k-1},
%\qquad
%E(R)=\frac{k-1}{\Lambda}(R_\infty-R).
\]
so that \eqref{eq:tildeKval} boils down to
\[
\widetilde K(R) \precsim  %\frac{k-1}{\Lambda}
R_\infty - R.
\]
This holds if we choose
\[
\widetilde K(R) = \frac{n\Lambda^2}{\lambda\gamma(R)} \quad \text{and} \quad\gamma(R)\succsim\frac1{R_\infty-R}.
\]

\noindent
{\bf Example.} Assume $k = 2$, $C/\Lambda\equiv c_{\Lambda}>0$, and
constant ellipticity bounds $0<\lambda\leq\Lambda$.
The associated ODE is
\[
f_\nu'=c_{\Lambda}f_\nu^2,
\qquad f_\nu(0)=\nu,
\]
and hence
\[
f_\nu(R)=\frac{\nu}{1-c_{\Lambda}\nu R},
\quad
h_\nu(R)=-\frac{1}{c_{\Lambda}}\log(1-c_{\Lambda}\nu R),
\quad
R_\infty=\frac{1}{c_{\Lambda}\nu}.
\]
Taking
\[
\widetilde K(R)=\frac{n\Lambda}{c_{\Lambda}\lambda\gamma(R)},
\quad \text{and} \quad
\gamma(R)\succsim \frac{1}{\frac{1}{c_{\Lambda}\nu}-R}
\]
gives
\[
\liminf_{R\rightarrow {\frac{1}{c_{\Lambda}\nu}}^-}
\frac{M(R)}
{\left|\log\left(\frac{1}{c_{\Lambda}\nu} - R\right)\right|} > 0.
\]
Thus, the growth is logarithmic and infinity is reached
at distance $R_\infty = 1/(c_{\Lambda}\nu)$ from the boundary.
We conclude sharpness by the function
\[
u(x)=-\frac{1}{c_{\Lambda}}\log(1-c_{\Lambda}\nu x_n)
\]
since it solves $\Delta u=c_{\Lambda}|Du|^2$ in the strip
$0<x_n<1/(c_{\Lambda}\nu)$ and vanishes on $x_n=0$, and satisfies
$M(R)=h_\nu(R)$.

\bigskip
\noindent
{\bf Example.} Assume again $C/\Lambda\equiv c_{\Lambda}>0$
and constant ellipticity bounds $0<\lambda\leq\Lambda$.
Straightforward calculations, as in the preceding example,
give the following growth scales depending on $k$:
\[
\begin{array}{c|c|c}
k & R_\infty & \text{Behaviour of }h_\nu(R)\\ \hline
k<1 & \infty
& h_\nu(R)\asymp
\displaystyle\frac{1}{c_{\Lambda}}
(c_{\Lambda}R)^{(2-k)/(1-k)}\\[6pt]
k=1 & \infty
& h_\nu(R)\asymp e^{c_{\Lambda}R}\\[3pt]
1<k<2 & \displaystyle\frac{\nu^{1-k}}{c_{\Lambda}(k-1)}
& h_\nu(R)\asymp
\displaystyle\frac{1}{c_{\Lambda}}
\bigl(c_{\Lambda}(R_\infty-R)\bigr)^{-(2-k)/(k-1)}\\[6pt]
k=2 & \displaystyle\frac{1}{c_{\Lambda}\nu}
& h_\nu(R)\asymp
\displaystyle\frac{1}{c_{\Lambda}}|\log(R_\infty-R)|\\[6pt]
k>2 & \displaystyle\frac{\nu^{1-k}}{c_{\Lambda}(k-1)}
& h_\nu(R)\to\displaystyle\frac{\nu^{2-k}}{c_{\Lambda}(k-2)}.
\end{array}
\]
Here $R\rightarrow R_\infty$.
Choosing
$\gamma(R)\succsim R$ for $k\leq1$ and
$\gamma(R)\succsim(R_\infty-R)^{-1}$ for $k>1$
gives a positive lower bound for $M'(R)/h_\nu'(R)$.
For $k\leq2$, integration yields the corresponding
lower growth bound for $M(R)$.
For $k>2$, the theorem gives
\[
\liminf_{R\rightarrow R_\infty}
M'(R)(R_\infty-R)^{1/(k-1)}>0,
\]
although $M(R)$ need not become unbounded.
These scales are attained by $u(x)=h_\nu(x_n)$,
which solves $\Delta u=c_{\Lambda}|Du|^k$ wherever defined.

\bigskip

\noindent\textbf{Example.}
We now consider $C = C(t)$ producing transitions
between linear growth, faster growth, and blow-ups.
Assume that $\lambda,\Lambda$
are constant and take
\[
\frac{C(t)}{\Lambda}=\frac{c_{\Lambda}}{(1+t)^{p}},
\quad c_{\Lambda}>0,\quad p>1,\quad k>1.
\]
This coefficient satisfies the structural assumptions with $t=x_n$.
The associated ODE gives
\[
f_\nu(t)=
\left[
\nu^{1-k}
-\frac{c_{\Lambda}(k-1)}{p-1}
\left(1 - \frac{1}{(1+t)^{
p-1}}\right)
\right]^{-1/(k-1)}.
\]
Introducing the critical initial slope
\[
\nu_*=
\left(\frac{p-1}{c_{\Lambda}(k-1)}\right)^{1/(k-1)},
\]
we obtain three distinct regimes:
\[
\begin{array}{c|c|c}
 & R_\infty & \text{Behaviour of }h_\nu(R)\\ \hline
0<\nu<\nu_* & \infty
& h_\nu(R)\asymp R\\[3pt]
\nu=\nu_* & \infty
& h_\nu(R)\asymp R^{(k + p - 2)/(k-1)}\\[3pt]
\nu>\nu_* & R_\infty<\infty
& \text{as in the constant-coefficient example}.
\end{array}
\]
In the last row, the same growth exponents and the
same classification according to $1<k<2$, $k=2$,
and $k>2$ apply, with $c_\Lambda$ replaced by
$c_\Lambda/(1+R_\infty)^p$ in the leading asymptotic
terms. For $k>2$, this describes the approach to the
finite limit, rather than the value of the limit itself.
The blow-up distance is
\[
R_\infty=
\left[1-\left(\frac{\nu_*}{\nu}\right)^{k-1}\right]^{-1/(p-1)}-1.
\]

At the critical slope $\nu_*$, the formulas simplify to
\[
f_{\nu_*}(t)=\nu_*(1+t)^{(p-1)/(k-1)},
\]
and
\[
h_{\nu_*}(R)=
\frac{(k-1)\nu_*}{k+p-2}
\left[(1+R)^{(k+p-2)/(k-1)}-1\right].
\]
Thus, the decay of $C$ shifts gradient blow-up to infinity
and produces an exact power growth rate.
For instance, if $p=k=2$, then $\nu_*=1/c_{\Lambda}$ and
\[
h_{\nu_*}(R)=\frac{1}{c_{\Lambda}}
\left(R+\frac{R^2}{2}\right).
\]
Concerning sharpness, we note that the function
\[
u(x)=\frac{1}{c_{\Lambda}}
\left(x_n+\frac{x_n^2}{2}\right)
\quad\text{solves}\quad
\Delta u=\frac{c_{\Lambda}|Du|^2}{(1+x_n)^2}.
\]

\bigskip

\noindent
\textbf{Example.}
We next consider an exponentially decaying coefficient,
producing another explicit critical regime.
Assume that $\lambda,\Lambda$ are constant and take
\[
\frac{C(t)}{\Lambda}=c_{\Lambda}e^{-at},
\quad c_{\Lambda}>0,\quad a>0,\quad k>1.
\]
This coefficient satisfies the structural assumptions
with $t=x_n$. The associated ODE gives
\[
f_\nu(t)=
\left[
\nu^{1-k}
-\frac{c_{\Lambda}(k-1)}{a}
\left(1-e^{-at}\right)
\right]^{-1/(k-1)}.
\]
Introducing the critical initial slope
\[
\nu_*=
\left(\frac{a}{c_{\Lambda}(k-1)}\right)^{1/(k-1)},
\]
we again obtain three distinct regimes:
\[
\begin{array}{c|c|c}
 & R_\infty & \text{Behaviour of }h_\nu(R)\\ \hline
0<\nu<\nu_* & \infty
& h_\nu(R)\asymp R\\[3pt]
\nu=\nu_* & \infty
& h_\nu(R)\asymp e^{aR/(k-1)}\\[3pt]
\nu>\nu_* & R_\infty<\infty
& \text{as in the constant-coefficient example}.
\end{array}
\]

In the last row, the same growth exponents and the
same classification according to $1<k<2$, $k=2$,
and $k>2$ apply, with $c_\Lambda$ replaced by
$c_\Lambda e^{-aR_\infty}$ in the leading asymptotic
terms. For $k>2$, this describes the approach to the
finite limit, rather than the value of the limit itself.
The blow-up distance is
\[
R_\infty=
-\frac1a\log\left[
1-\left(\frac{\nu_*}{\nu}\right)^{k-1}
\right].
\]

At the critical slope $\nu_*$, the formulas simplify to
\[
f_{\nu_*}(t)=\nu_*e^{at/(k-1)}
\]
and
\[
h_{\nu_*}(R)=\frac{(k-1)\nu_*}{a}
\left(e^{aR/(k-1)}-1\right).
\]
Thus, the decay of $C$ shifts gradient blow-up to infinity
and produces an exact exponential growth rate.
For instance, if $a=1$ and $k=2$, then
$\nu_*=1/c_{\Lambda}$ and
\[
h_{\nu_*}(R)=\frac1{c_{\Lambda}}(e^R-1).
\]
Concerning sharpness, we note that the function
\[
u(x)=\frac1{c_{\Lambda}}(e^{x_n}-1)
\quad\text{solves}\quad
\Delta u=c_{\Lambda}e^{-x_n}|Du|^2.
\]

In each case above, choosing $\gamma = \gamma(R)$ according to the preceding
criteria yields the corresponding lower growth estimates
from Theorem \ref{th:Phragmen--Lindelof}.
The estimates are attained by $u(x)=h_\nu(x_n)$, which solves
$\Lambda\Delta u=C(x_n)|Du|^k$ wherever defined.

%%%%%%%%%%%%%%%%%%%%%%%%%%%%%%%%%%%%%%%%%%%%
%%%%%%%%%%%%%%%%%%%%%%%%%%%%%%%%%%%%%%%%%%%%
%%%%%%%%%%%%%%%%%%%%%%%%%%%%%%%%%%%%%%%%%%%%

\subsection{The case $\Phi(|x|,u,|Du|)=-c-bu-a|Du| \leq 0$}
\label{subsec:dependonsolution}

We assume that $a,b,c\geq0$ are constants and first consider
$a+b>0$.
The ODE \eqref{eq:ODE} under consideration now reads
\begin{align*}
\frac{d^2h}{dt^2}
= \frac{c+b h(t)+a h'(t)}{\Lambda(t)}
-K(R)\frac{\Lambda(t)}{\lambda(t)}h'(t),
\quad t>0,\quad h(0)=0,\quad h'(0)=\nu>0.
\end{align*}
In the classical case of $\lambda$ and $\Lambda$ constant,
this reduces to the following ODEs for $h_{\nu,R}$ and $h_\nu$:
\[
h_{\nu,R}''(t)
-\left(\frac{a}{\Lambda}
-K(R)\frac{\Lambda}{\lambda}\right)h_{\nu,R}'(t)
-\frac{b}{\Lambda}h_{\nu,R}(t)
=\frac{c}{\Lambda},
\]
and
\[
h_\nu''(t)-\frac{a}{\Lambda}h_\nu'(t)
-\frac{b}{\Lambda}h_\nu(t)=\frac{c}{\Lambda}.
\]
Standard calculus then yields, for sufficiently large $R$,
\[
h_{\nu,R}(t)=
\frac{1}{\bar r_+-\bar r_-}
\left[
\nu\left(e^{\bar r_+t}-e^{\bar r_-t}\right)
+\frac{c}{\Lambda}
\left(
\frac{e^{\bar r_+t}-1}{\bar r_+}
-\frac{e^{\bar r_-t}-1}{\bar r_-}
\right)
\right]
\]
and
\[
h_\nu(t)=
\frac{1}{r_+-r_-}
\left[
\nu\left(e^{r_+t}-e^{r_-t}\right)
+\frac{c}{\Lambda}
\left(
\frac{e^{r_+t}-1}{r_+}
-\frac{e^{r_-t}-1}{r_-}
\right)
\right],
\]
where
\[
\bar r_\pm=
\frac{
a-K\frac{\Lambda^2}{\lambda}
\pm\sqrt{\left(a-K\frac{\Lambda^2}{\lambda}\right)^2+4b\Lambda}
}{2\Lambda}
\quad\text{and}\quad
r_\pm=\frac{a\pm\sqrt{a^2+4b\Lambda}}{2\Lambda}.
\]
Here $(e^{rt}-1)/r$ is interpreted as $t$ when $r=0$.
Since $a,b\geq0$ and $h_{\nu,R}'>0$, comparison of
the two linear ODEs gives $h_{\nu,R}<h_\nu$ for $t>0$.

As $r_+>0\geq r_-$ and $\bar r_\pm\to r_\pm$,
we have, for large enough $R$,
\begin{align*}
\frac{h'_{\nu,R}(R)}{h'_\nu(R)}
&=
\frac{
(\nu\bar r_++c/\Lambda)e^{\bar r_+R}
-(\nu\bar r_-+c/\Lambda)e^{\bar r_-R}
}{
(\nu r_++c/\Lambda)e^{r_+R}
-(\nu r_-+c/\Lambda)e^{r_-R}
}
\frac{r_+-r_-}{\bar r_+-\bar r_-}\\
&=(1+o(1))e^{(\bar r_+-r_+)R}
\geq\frac14e^{(\bar r_+-r_+)R}.
\end{align*}
Hence,
\[
\liminf_{R\to\infty}
\frac{h'_{\nu,R}(R)}{h'_\nu(R)}
\]
is positive if $(\bar r_+-r_+)R$ is bounded below.
This holds for the choice $\gamma(R)=R$, giving
$K=K(R)=n/R$.
Indeed, since $x\mapsto\sqrt{x^2+4b\Lambda}$
is $1$-Lipschitz, we have
\begin{align*}
0\geq(\bar r_+-r_+)R
&=\frac{R}{2\Lambda}
\left(
\sqrt{\left(a-K\frac{\Lambda^2}{\lambda}\right)^2+4b\Lambda}
-\sqrt{a^2+4b\Lambda}
-K\frac{\Lambda^2}{\lambda}
\right)\\
&\geq\frac{R}{2\Lambda}
\left(-K\frac{\Lambda^2}{\lambda}
-K\frac{\Lambda^2}{\lambda}\right)
=-n\frac{\Lambda}{\lambda}>-\infty.
\end{align*}
Thus, we have proved
\[
\liminf_{R\to\infty}
\frac{M'(R)}{h'_\nu(R)}>0
\quad\text{and hence}\quad
\liminf_{R\to\infty}
\frac{M'(R)}{e^{r_+R}}>0.
\]
Integration shows that the subsolution to the PDE
must grow at least as fast as $e^{r_+R}$, i.e.,
\[
\liminf_{R\to\infty}\frac{M(R)}{e^{r_+R}}>0.
\]
The estimate is sharp, since $u(x)=h_\nu(x_n)$ solves
\[
\Lambda\Delta u=c+bu+a|Du|
\]
in $\mathbb R^n_+$, vanishes on $x_n=0$,
and satisfies $M(R)=h_\nu(R)\asymp e^{r_+R}$.

Finally, if $a=b=0$, then
\[
h_\nu(t)=\nu t+\frac{c}{2\Lambda}t^2.
\]
For $c>0$, the same choice $\gamma(R)=R$ gives
\[
\lim_{R\to\infty}
\frac{h'_{\nu,R}(R)}{h'_\nu(R)}
=
\frac{1-e^{-n\Lambda/\lambda}}{n\Lambda/\lambda}>0,
\]
and hence the sharp quadratic growth estimate
\[
\liminf_{R\to\infty}\frac{M(R)}{R^2}>0.
\]
If also $c=0$, the corresponding sharp growth is linear.

%%%%%%%%%%%%%%%%%%%%%%%%%%%%%%%%%%%%%%%%%%%%%%%%%%%%%%%%%%%%%%%%%%%%%
%%%%%%%%%%%%%%%%%%%%%%%%%%%%%%%%%%%%%%%%%%%%%%%%%%%%%%%%%%%%%%%%%%%%%
%%%%%%%%%%%%%%%%%%%%%%%%%%%%%%%%%%%%%%%%%%%%%%%%%%%%%%%%%%%%%%%%%%%%%
%%%%%%%%%%%%%%%%%%%%%%%%%%%%%%%%%%%%%%%%%%%%%%%%%%%%%%%%%%%%%%%%%%%%%
%%%%%%%%%%%%%%%%%%%%%%%%%%%%%%%%%%%%%%%%%%%%%%%%%%%%%%%%%%%%%%%%%%%%%
%%%%%%%%%%%%%%%%%%%%%%%%%%%%%%%%%%%%%%%%%%%%%%%%%%%%%%%%%%%%%%%%%%%%%
%%%%%%%%%%%%%%%%%%%%%%%%%%%%%%%%%%%%%%%%%%%%%%%%%%%%%%%%%%%%%%%%%%%%%
%%%%%%%%%%%%%%%%%%%%%%%%%%%%%%%%%%%%%%%%%%%%%%%%%%%%%%%%%%%%%%%%%%%%%
%%%%%%%%%%%%%%%%%%%%%%%%%%%%%%%%%%%%%%%%%%%%%%%%%%%%%%%%%%%%%%%%%%%%%
%%%%%%%%%%%%%%%%%%%%%%%%%%%%%%%%%%%%%%%%%%%%%%%%%%%%%%%%%%%%%%%%%%%%%
%%%%%%%%%%%%%%%%%%%%%%%%%%%%%%%%%%%%%%%%%%%%%%%%%%%%%%%%%%%%%%%%%%%%%
%%%%%%%%%%%%%%%%%%%%%%%%%%%%%%%%%%%%%%%%%%%%%%%%%%%%%%%%%%%%%%%%%%%%%
%%%%%%%%%%%%%%%%%%%%%%%%%%%%%%%%%%%%%%%%%%%%%%%%%%%%%%%%%%%%%%%%%%%%%
%%%%%%%%%%%%%%%%%%%%%%%%%%%%%%%%%%%%%%%%%%%%%%%%%%%%%%%%%%%%%%%%%%%%%
%%%%%%%%%%%%%%%%%%%%%%%%%%%%%%%%%%%%%%%%%%%%%%%%%%%%%%%%%%%%%%%%%%%%%

\end{document}